\documentclass[12pt,psamsfonts]{article}

\usepackage{theorem}
\usepackage{amscd}
\usepackage{babel}
\usepackage[leqno]{amsmath}
\usepackage{amssymb}
\usepackage{amsfonts}

\theoremstyle{change}
\theorembodyfont{\normalfont}
\newtheorem{atom}       {}      [subsection]
\newtheorem{definition} [atom]{Definition.}
\newtheorem{remark}     [atom]{Remark. }

\newtheorem{example}    [atom]{Example.  }

\theoremheaderfont{\normalfont\bfseries\rmfamily}
\theorembodyfont{\normalfont\itshape}
\newtheorem{lemma}        [atom]{Lemma.}
\newtheorem{proposition}  [atom]{Proposition.} 
\newtheorem{theorem}      [atom]{Theorem.}
\newtheorem{corollary}    [atom]{Corollary.}

\newcommand{\ZZ}{\mathbb{Z}}
\newcommand{\QQ}{\mathbb{Q}}
\newcommand{\RR}{\mathbb{R}}

\newcommand{\PP}{\mathbb{P}}

\newcommand  {\shE}     {\mathcal{E}}

\newcommand  {\shI}     {\mathcal{I}}

\newcommand  {\shM}     {\mathcal{M}}

\newcommand  {\shN}     {\mathcal{N}}
\newcommand  {\shL}     {\mathcal{L}}

\newcommand  {\foX}     {\mathfrak{X}}

\newcommand  {\aff}     {{\text{aff}}}

\newcommand  {\cyc}     {\operatorname{cyc}}

\renewcommand{\div}     {\operatorname{div}}
\newcommand  {\Div}     {\operatorname{Div}}
\newcommand  {\dom}     {\operatorname{dom}}

\newcommand  {\Ext}     {\operatorname{Ext}}

\renewcommand  {\k}       {\kappa}

\newcommand  {\lra}     {\longrightarrow}

\newcommand  {\NAP}     {\operatorname{\overline{NA}}}

\newcommand  {\NEP}     {\operatorname{\overline{NE}}}

\renewcommand{\O}       {\mathcal{O}}

\newcommand  {\Pic}     {\operatorname{Pic}}

\newcommand  {\proof}   {\textit{Proof. }}
\newcommand  {\Proj}    {\operatorname{Proj}}

\newcommand  {\ra}      {\rightarrow}
\newcommand  {\Ra}      {\Rightarrow}

\newcommand  {\Reg}     {\operatorname{Reg}}

\newcommand  {\Sing}    {\operatorname{Sing}}

\newcommand  {\Spec}    {\operatorname{Spec}}
\newcommand  {\Supp}    {\operatorname{Supp}}
\newcommand  {\Sym}     {\operatorname{Sym}}

\begin{document}

\setcounter{subsection}{0}
\setlength{\unitlength}{1ex}

\begin{titlepage}
\title {On contractible curves on   normal surfaces}
\author{Stefan Schr\"oer\\
        }
\end   {titlepage}


\maketitle

\begin{abstract}
We give  characterizations of contractible curves
on  proper normal algebraic surfaces in terms of complementary
Weil divisors. From this we obtain some generalizations of
the classical criteria for contractibility of Castelnuovo and
Artin. Furthermore, we will derive a finiteness result on
homogeneous spectra defined by Weil divisors on proper normal
algebraic surfaces.
\end{abstract}


\renewcommand{\thefootnote}{}
\footnote{\hspace{-4ex}{\itshape Key words:} normal surfaces,
contraction, modification.}
\footnote{\hspace{-4ex}{\itshape Mathematics subject 
classification (1991)}:   14C20, 14E10, 14J05, 14J17.}

\subsection{Introduction}

The goal of this paper is to characterize contractible curves
on proper normal  algebraic surfaces.

The corresponding question for compact normal complex-analytic
surfaces was completely solved by Grauert
\cite{Grauert 1962}: A curve 
$R $ on such a complex-analytic surface 
$X $ is contractible
to another   complex-analytic surface 
$Y $ if and only if the curve
$R $ is \textit{negative definite}.
Later, Artin
\cite{Artin 1970} showed that a similar result holds in the
category of proper normal 2-dimensional algebraic spaces.

Usually,  if one  starts with a  proper normal
algebraic surface 
$X $ and a negative definite curve 
$R\subset X $, the resulting algebraic space 
$Y $ is not a scheme anymore. We want to characterize those
curves 
$R\subset X $ for which 
$Y $ remains a scheme. For simplicity, we call such curves
\textit{contractible}.  Our question was already  considered by
Artin in a series of two papers 
\cite{Artin 1962} 
\cite{Artin 1966}, where he gives several necessary or
sufficient conditions for smooth surfaces. For example, the
classical Castel\-nuovo criterion  tells us that an exceptional
curve of the first kind is contractible.

\medskip
Our main result is  a characterization of contractible curves 
$R\subset X $ in
terms of complementary Weil divisors, see theorem
\ref{complementary cycles}. Roughly speaking, one has to find
a Weil divisor which defines the contraction numerically  and
behaves reasonably on the formal completion 
$\foX $ of 
$R\subset X $. From this
characterization we will easily derive  some generalizations of
the classical results of Castel\-nuovo and Artin. Furthermore, our
approach will explain certain  differences between characteristic 
$p>0 $ and characteristic zero.

The results will be useful in the classification \`a la
Enriques of proper normal algebraic surfaces, since they allow us
to introduce the technique of extremal rays from Mori theory
without conditions on the singularities. Since our ground field 
$K $ is arbitrary, the results might also be applied   to generic
hyperplane sections or generic fibres in higher dimensional
geometry.

As an application we will show that for every Weil divisor 
$D\in Z^1(X) $, the homogeneous spectrum  
$P(X,D) $ of the graded algebra 
$  \oplus_{n\geq 0}H^0(X,nD) $ is always a scheme of finite
type. Here the interesting point is that the \textit{algebra} 
usually fails to be finitely generated.

\medskip
This article is divided into six sections, including the
introduction.   In the second section we  will set up the
notation and recall some facts about the cone of curves for
proper  normal surfaces. The third section contains a
characterization of almost affine open subsets and  our main
results on contractible curves. The fourth section is a bit
technical. We introduce the notion of   numerical 
$\QQ $-factoriality and look for  situations in which   our
criteria are applicable. The fifth section contains some
sufficient conditions for contractibility. In the last section we
prove that the homogeneous spectra
$P(X,D) $ are always  schemes of finite type.

\medskip
I would like to thank the referee for pointing out an error 
in a preliminary version and  several useful suggestions.

\subsection{The cone of curves}

This section is preparatory in nature. We will set up the
notation and discuss the   cone of curves for proper normal
surfaces. Throughout the paper, we will work over an arbitrary
base field 
$K $. The  word \textit{surface}  always
refers to a 2-dimensional, irreducible, separated
$K $-scheme of finite type. In this section, 
$X $ will always be a proper normal surface.

\begin{atom}
Let 
$Z^1(X) $ be the group of Weil divisors.
If 
$X $ is a proper regular surface, the group of Weil divisors comes
along with the 
$\ZZ $-valued intersection pairing 
$A,B\mapsto A\cdot B $. If 
$X $ is just a proper normal surface, Mumford
\cite{Mumford  1961}
pointed out that there is still a 
$\QQ $-valued intersection form on 
$Z^1(X) $ with the usual properties. Mumford's approach uses
resolutions of singularities for surfaces, which has been
established in full generality by Lipman
\cite{Lipman 1978}. 
\end{atom}

\begin{atom}
\label{negative definite curves}
We will call a curve 
$R\subset X $ \textit{negative definite}  if the intersection
matrix 
$\Phi=(R_i\cdot R_j) $ for the irreducible components 
$R_i\subset R $ is negative definite. 
For example, if
$f:X\ra Y $ is  a birational morphism of proper normal surfaces,
the contracted curve 
$R\subset X $ is negative definite, see for example
\cite{SGA 7}, X, 1.9. Since  
$R_i\cdot R_j\geq 0 $ holds for 
$i\neq j $, it is a well known fact from linear algebra that 
the inverse matrix
$\Phi^{-1} $ has negative entries. Moreover, if  
$R $ is connected, all entries are strictly negative. We will
use this fact several times. 
\end{atom}

\begin{atom}
A Weil divisor 
$A\in Z^1(X) $ which lies in the radical of the rational  
intersection pairing is called \textit{numerically trivial}.  Let
us write 
$ N(X) $ for the quotient of 
$Z^1(X) $ modulo  the radical of  Mumford's
rational intersection pairing.
We obtain a
non-degenerate intersection pairing 
$$
N(X)\times  N(X)\lra \QQ. 
$$
It is well-known that 
$ N(X) $ is a free 
$\ZZ $-module of finite type. The Hodge
index theorem tells us that this  intersection form has exactly
one positive eigenvalue. 
In Mori theory, one traditionally works with the pairing 
$N^1(X)\times N_1(X)\ra \ZZ $ induced
by
$\Div(X)\times Z^1(X)\ra \ZZ $. But 
$\Pic(X)=0$ might easily happen if 
$X $ is non-projective
\cite {Schroeer b}, in which case this pairing has lost its
significance.
\end{atom}
 
\begin{atom}
We call a subset of a real vector space a {\em cone} if it  is
closed under multiplication with positive scalars and
addition. The real vector space 
$N(X,\RR)=N(X)\otimes\RR $ contains two important cones,
which we should recall.

The closed cone 
$\NEP(X) $ generated by all curves 
$C\subset X$ is called the
{\em pseudo-effective} cone.
 
A Weil divisor 
$A $ meeting the condition of the Nakai criterion for ampleness,
which means 
$A^2>0 $ and 
$ A\cdot C>0$ for all curves 
$C\subset X $, will be called
{\itshape ample}. The {\em pseudo-ample}   cone 
$\NAP(X) $ is the closed cone generated by  the ample cycles.
Note that an ample Weil divisor  does not define a closed
embedding into some 
$\PP^n$ unless it is a 
$\QQ $-Cartier divisor. 

It follows from
\cite{SGA 6}, XIII, 7.1 that  
$\NEP(X) $ is polar to 
$\NAP(X)$. 
\end{atom}

\medskip
Since the pseudo-effective cone is closed and contains no lines, 
it is generated by its extremal rays. One calls a closed subcone 
$P\subset \NEP(X) $ {\itshape extremal} if 
$e+e'\in P $ for 
$e,e'\in \NEP(X) $ implies 
$e,e'\in P $. The extremal closed subcones  have the following
geometric significance:

\begin{proposition}
\label{extremal cones}
The extremal closed subcones 
$P\subset \NEP(X) $ are precisely the following subsets:
\renewcommand{\labelenumi}{(\roman{enumi})}
\begin{enumerate}
\item Cones of the form
$P=\sum\RR_+ R_i $, generated  by the irreducible components 
$R_i $ of
a unique reduced negative definite curve 
$R\subset X $.
\item Cones of the form 
$P=\NEP(X)\cap F^\perp $ for certain real pseudo-effective
class 
$F\in \partial \NAP(X) $ with 
$F^2=0 $.
\end{enumerate}
  Moreover, the cones 
$P=\sum\RR_+ R_i $ of the first sort satisfy the following 
finiteness condition: If 
$Q\subset N(X,\RR) $ is the closed cone generated by all 
irreducible curves 
$ C\subset X$ not supported by 
$R $, then 
$P\cap Q=0 $ holds.
\end{proposition}

\proof 
An extremal closed subcone 
$P\subset \NEP(X) $  
is of the form $ \NEP(X) \cap a^\perp$ for some support
function 
$a\in \partial\NAP(X) $. If 
$ a^2=0$ holds, we are in case 
$(ii) $. Assume that 
$a^2>0 $ holds. Then the intersection form on the vector
subspace
$P-P $ is negative definite. Write 
$P=\sum P_i$ as the convex hull of   extremal rays.
The argument in 
\cite{Kollar 1995},  lemma 4.12 for smooth surfaces also
applies for normal surfaces, consequently each 
$P_i $ is generated by a unique integral, negative definite
curve 
$R_i\subset X $. 
If there would be a relation 
$e_j=\sum_{i\neq j}\lambda_ie_i $, the coefficients must be
positive, and we obtain the contradiction 
$$
0>R_j^2 = \sum_{i\neq j}\lambda_i R_i\cdot R_j \geq 0. 
$$
So the  
$R_i $ are linearly independent, hence finite in number, and 
$R=\cup R_i $ is the negative definite curve whose irreducible
components generate 
$P $.
Conversely, one easily sees that the irreducible components of a
negative definite curve generate an extremal cone. 

To verify the finiteness condition, let 
$e\in P\cap Q $. We have 
$e=\sum \lambda_i R_i$ with positive coefficients  and 
$e=\lim C_n $ for certain  real effective 1-cycles 
$C_n $ whose components are not contained in 
$R $. Hence we have 
$$
0>e^2 =\lim(e\cdot C_n)  = \lim(\sum\lambda_iR_i\cdot C_n) \geq 0, 
$$
contradiction. QED.

\begin{remark}
It is possible to extend Mumford's intersection pairing
fruitfully to {\itshape unibranched} surfaces. Let 
$X $ be a unibranched surface, which means that the normalization
$ X'\ra X $ is a bijective map
\cite{EGA IV 1}, EGA 0, 23.2.1. Obviously  we have 
$Z^1(X)=Z^1(X') $. If 
$\eta\in X $ is the generic point, let 
$d $ be the length of the Artin ring 
$\O_{X,\eta} $. For two Weil divisors 
$A,B $ on 
$X $ the result 
\cite{EGA IV 4}, EGA IV, 21.10.4 forces us to put 
$$
A\cdot B = d( A\cdot B)_{X'}, 
$$
where the right hand side is computed  on the normalization.
The results of this paper remain true for unibranched
surfaces, but we are content with the normal case. 
Nevertheless,  it should be noted that from a conceptual point
of view, unibranched surfaces form  a better category, since this
category is stable under extensions of the ground field
$K$, which is not true for normal surfaces.
\end{remark}

\subsection{Characterization of contractible curves}

This section contains our main results on contractible curves.
Throughout,
$X $ will be a proper normal surface. We first make the trivial
observation that it suffices to treat the case of {\itshape
connected}  curves:

\begin{lemma} 
\label{connected   curves}
A curve 
$R\subset X $ is contractible if and only if all  its connected
components 
$R_i\subset X $ are contractible.
\end{lemma} 
 
\proof Use patching. QED.

\medskip
For the purpose of this paper, the following terminology will
be useful.   We call an  open
subset 
$U\subset X $ {\em almost affine}  if the affine hull
$U^\aff=\Spec\, \Gamma(U,\O_{X}) $ is of finite type over  the
ground field 
$K $ and the canonical morphism
$U\ra U^\aff $ is proper and birational. We remark that these are
precisely the semi-affine open subsets of Goodman and Landman
\cite{Goodman-Landman 1973} with a 2-dimensional ring of global
sections
$\Gamma(U,\O_{X}) $. There is a kind of dual
notion for curves: We say that a   curve 
$A\subset X $ is {\em ample on itself} if 
$A\cdot A_i>0 $ holds for all irreducible components 
$A_i\subset A $. These concepts are related in the following
way:

\begin{proposition} 
\label{almost affine open subset}
Let 
$ U\subset X$ be an open subset, 
$C=X\setminus U $ its complement, and 
$C_i\subset C $ the integral components. Then the following
conditions are equivalent:
\renewcommand{\labelenumi}{\alph{enumi})}
\begin{enumerate}
\item 
The open subset 
$U \subset X$ is almost affine.
\item  
The closed subset 
$C\subset X $ is a connected curve and some Weil divisor 
$A$ supported by 
$C$ has 
$A^2>0$.
\item 
There is a   curve 
$A\subset X $ with  
$\Supp(A)= C$ which is ample on itself.
\end{enumerate}
\end{proposition}

\proof
We verify the implications 
$(a)\Ra (b)\Ra (c)\Ra   (a)$. 
Assume that 
$ U$ is almost affine. Passing to the contraction 
$f:X\ra Y $ defined by 
$U\ra U^\aff $, we can assume that 
$U $ is affine. Then 
$\O_U $ is ample, and we easily find a section 
$s\in\Gamma(U,\O_{X}) $ defining a Cartier divisor 
$D_U\subset U $ whose closure 
$\overline{D_U}\subset X $ intersects each component 
$C_i $. According to 
\cite{SGA 6}, II, 2.2.6, each 
$C_i $ is a curve. Multiplying 
$s $ by a suitable power  of   
$t\in \Gamma(U,\O_{X}) $ defined by the canonical homomorphism
$\O_{X}\ra\O_{X}(C) $ we can assume that the domain of definition 
$\dom(s) \subset X$ equals 
$U $. Let 
$D=\div(s) $ be the corresponding    principal  divisor on 
$X $ with 
$D\cap U=D_U $ and write 
$\cyc(D)=D_1-D_2 $ with 
$D_1=\overline{D_U} $; then 
$D_2\cdot C_i=D_1\cdot C_i>0 $ holds for each component 
$C_i $. 
Since 
$\dom(s)=U $, the Weil divisor 
$D_2 $ is effective with 
$\Supp(D_2)=C $. Suppose 
$C'\subset C $ is a connected component. Let
$D_2'\subset D_2 $ be the corresponding connected component.
Then 
$(D_2')^2=D_2'\cdot D_2 = D_2'\cdot D_1 >0 $. 
Moreover, the curve 
$C $ is be connected by the Hodge index theorem.
Hence $A=D_2$ satisfies condition
$(b)$.

Now assume that $(b)$ holds. Decompose 
$A=A_+ - A_-$ into positive and negative part. Then 
$0<A^2=A_+^2-2A_+A_- +A_-^2\leq A_+^2+A_-^2$.
So there is a curve 
$A_0\subset X $ supported by 
$C $ with 
$A_0^2>0 $.
Now we   rename the irreducible
components of $C$ and find a complete list 
$C_1,\ldots,C_n $, possibly with repetitions, with 
$A_0\cdot C_1>0 $ and
$C_i\cap C_{i+1}  $  nonempty. Inductively we define  curves 
$A_j\subset X $ with support 
$A_0\cup C_1\cup\ldots\cup C_j $ and 
$A_j\cdot C_i >0$ for 
$1\leq i\leq j $ as follows:  If 
$A_j $ is already defined, then 
$A_{j+1}=\lambda_jA_j+C_{j+1} $ meets our conditions for 
$\lambda_j>0 $ sufficiently large. Now 
$A=A_n $ is the desired   curve which is ample on itself.

Finally, assume that 
$ (c)$ holds. Let 
$g:X'\ra X $ be a resolution of all the singularities 
$x\in \Sing(X)\cap C $. Then 
$ g^*(A) $ is an effective
$\QQ $-Cartier divisor  with support 
$A'=g^{-1}(A) $. Let 
$E'\in\Div(X')$ be  a relatively ample divisor whose support is
contracted by 
$g $. Then 
$A'=g^*(nA)+E' $ is effective with support 
$A' $, provided 
$n>0 $ is sufficiently large, and ample on itself. Hence it
suffices to treat the case that 
$A $ is a Cartier divisor. The associated invertible sheaf 
$\shL=\O_{X}(A) $ is ample on 
$A $. According to   
\cite{Fujita 1983}, theorem 1.10, a multiple 
$\shL^{\otimes n} $ is globally generated, and the homogeneous
spectrum $Y$ of 
$\Gamma(X,\Sym\shL) $ yields a birational contraction 
$ f:X\ra Y$. We infer that 
$U=f^{-1}(V) $ is the preimage of an affine open subset 
$V\subset Y $, hence is almost affine. QED.

\medskip
From this we easily derive a geometric characterization of
contractible curves in terms of complementary  
curves ample on themselves:

\begin{theorem}
\label{complementary curves}
A connected negative definite curve 
$R\subset X $ is contractible if and only if there is a curve 
$A\subset X $ disjoint to 
$R $ with 
$A\cdot C>0 $ for every curve 
$C\subset X $ not supported by 
$R $.
\end{theorem}
 
\proof
Assume there is a contraction 
$f:X\ra Y $ mapping 
$R $ to a closed point 
$ y\in Y$. Then the complement 
$ Y\setminus V $ of an affine open neighborhood 
$V\subset Y$ can be regarded as a curve on 
$X $. According to 
\ref{almost affine open subset} it is the  support of  a curve 
$A\subset X $ ample on itself with the desired property.
Conversely, if there is such a curve 
$A\subset X$, the open subset 
$U=X\setminus A $ is almost affine and defines the contraction
$f:X\ra Y $ of
$R $. QED.

\medskip
We have seen that one can always find a possibly non-effective
Weil divisor 
$A $ satisfying the above condition, and the real problem is to
choose an {\em effective} one. Our main result has the
advantage that it does not involve any conditions of
effectivity:

\begin{theorem}
\label{complementary cycles}
A connected negative definite curve
$R\subset X $   is contractible if and only if there is  a
Weil divisor
$A\in Z^1(X) $ satisfying the following three conditions:

\renewcommand{\labelenumi}{(\roman{enumi})}
\begin{enumerate}
\item The Weil divisor 
$A $ is Cartier  near 
$R\subset X $.
\item We have 
$A\cdot C\geq 0 $ for all curves 
$C\subset X $, with equality if and only if 
$C\subset R $ holds.
\item For every integer 
$m>0 $ there is an integer 
$ n>0$ and a numerically trivial Weil divisor 
$ N\in Z^1(X)$ which is Cartier near 
$ R\subset X$ such that the linear class of 
$ nA+N $ is trivial on 
$mR$.
\end{enumerate}
\end{theorem}
 
\proof
According to 
\ref{complementary curves}
the conditions are necessary, and the real issue is to show
sufficiency.

We first get rid of the singularities. Let 
$g:X'\ra X $ be a resolution of all singularities
$x\in \Sing(X)\cap R$. According to 
\cite{EGA II}, EGA II, 8.11.1 we have to show
that 
$R'=g^{-1}(R) $ is contractible.  A straightforward argument
shows that 
$A'=g^*(A) $ satisfies the three conditions on 
$X' $, hence we can assume that 
$R\subset \Reg(X) $ holds. Now let
$g:X'\ra X $ be a resolution of all remaining singularities.
According to 
\ref{connected curves}
we have to show that 
$R $ viewed as a curve 
$R'\subset X' $ is contractible. Let 
$P' \subset N(X',\RR) $ be the  cone generated by  the
irreducible components of
$R'$ and 
$Q'\subset N(X',\RR) $ the closed cone generated by all other 
integral curves. Choose a 
 relatively ample exceptional divisor 
  $D'\in Z^1(X' )$ contracted by $g$. If 
$U'\subset N(X',\RR) $ is a  compact neighborhood of zero, the
linear form associated to 
$D' $ is bounded on 
$U'\cap Q' $. Since 
$P'\cap Q'=0 $ holds by
\ref{extremal cones}, the 
$\QQ $-divisor 
$A'= D' + g^*(nA)$ is strictly positive on the punctured
neighborhood 
$U'\cap Q'\setminus\left\{  0  \right\}$ provided 
$n $ is sufficiently large. Consequently, the
numerical class of 
$A' $ is a support function of the pseudo-effective cone with
respect to the extremal subcone
$P' $, and condition $(i)$ holds true on 
$X' $. The other conditions trivially remain unaffected, hence
we can assume that 
$X $ is a  {\em regular}  surface.

Obviously 
$A^2>0 $ holds.
Since 
$R $ is negative definite, there exists  a divisor 
$D\in Z^1(X) $ supported by $R$ which  is anti-ample on 
$R $; according to 
\ref{negative definite curves}, the divisor is effective and its
support equals
$R $. We have 
$$
(tA-D)^2 = t^2A^2 - 2t A\cdot D + D^2 >0
$$
for all 
$t>0 $ sufficiently large. Arguing as above we also see that 
$(tA-D)\cdot C>0 $ holds for all   curves 
$C\subset Y $, provided 
$t>0 $ is sufficiently large.
We now replace 
$A$ by a suitable multiple and assume, using the Nakai
criterion, that 
$A-D $ is ample.

Set  
$\shL=\O_{X}(A) $ and 
$\shI=\O_{X}(-D) $; thus  
$\shL $ is pseudo-ample and 
$\shL\otimes\shI $ is ample. According to Fujita's  vanishing
result
\cite{Fujita 1983}, theorem 5.1, there is a natural number 
$t_0>0 $ such that 
$$
H^1(X,
\shL^{\otimes  s}\otimes (\shL\otimes \shI)^{\otimes t}\otimes
\shN) =0 
$$
holds for all integers 
$s\geq 0 $,
$t\geq t_0 $ and all numerically trivial invertible sheaves 
$\shN $.
If we replace 
$\shL $ and 
$\shI $ by 
$ \shL^{\otimes t_0}$ and 
$ \shI^{\otimes t_0}$ we can assume that 
$$
H^1(X,\shL^{\otimes s} \otimes \shI \otimes \shN)=0
$$
holds for all integers 
$s>0 $ and all numerically trivial invertible sheaves 
$\shN $.
Thus the right hand term in the exact sequence 
$$
H^0(X,\shL^{\otimes s}\otimes \shN) \lra
H^0(D,\shL^{\otimes s}\otimes \shN\mid D) \lra
H^1(X, \shL^{\otimes s}\otimes \shI\otimes\shN)
$$
is zero. Now choose an integer 
$m $ with  
$D\subset mR  $. According to  condition $(c)$, we can pick a
numerically trivial invertible sheaf 
$\shN $ and an integer 
$n>0 $ such that the restriction
$\shL^{\otimes n}\otimes \shN\mid D $ is  trivial. Let us
replace 
$\shL $ by 
$\shL^{\otimes n}\otimes \shN $. 
Consider the open subset 
$U\subset X $ on which 
$\shL $ is globally generated; by construction, this open set
contains the generic points of
$R $, hence 
$\shL $ is ample on the complement 
$X\setminus U $. According to 
\cite{Fujita 1983}, theorem 1.10, some multiple  
$\shL^{\otimes t} $ is globally generated, and the homogeneous
spectrum 
$$
Y= \Proj( \Gamma(X,\Sym\shL))
$$ 
is a projective normal surface yielding the
desired contraction
$f:X\ra Y$ of 
$R $. QED.

\subsection{Improvement of cycles}

How can one ensure that the conditions of theorem 
\ref{complementary cycles} are fulfilled? It is easy to find 
  a Weil divisor 
$A\in Z^1(X) $ satisfying the condition (ii) of 
\ref{complementary cycles}.  In this section we discuss the
possibilities to improve the Weil divisor 
$A $ in such a way that it also   satisfies the other two more
delicate  conditions (i) and (iii). For this purpose the following
notion is useful:

\begin{definition}
Let 
$S $ be an arbitrary subset of a proper normal surface 
$X $. We say that 
$X $ is {\em numerically 
$\QQ $-factorial} with respect to 
$S $ if   each Weil divisor
$D $ is numerically equivalent to a 
$\QQ $-Weil divisor which is 
$\QQ $-Cartier near 
$S $.
\end{definition}

Of course, only the points 
$x\in S $ whose local rings 
$\O_{X,x} $ are not 
$\QQ $-factorial are relevant for this notion. If the  conditions
holds for 
$S=X $  we call the surface
$X $ \textit{numerically 
$\QQ $-factorial}.  According to the Nakai criterion, such a 
surface
is projective; more precisely, the canonical map 
$\Pic(X)\ra N(X) $ has finite cokernel.

We have the following behaviour under birational morphisms:

\begin{proposition}
\label{factoriality and birational}
Let 
$f:X\ra Y $ be a birational morphism of proper normal surfaces, 
$R\subset X $ the contracted curve, 
$\foX \subset X$ the corresponding formal completion, 
$S\subset X $ a subset containing 
$R $, and 
$T\subset Y $ its image. If the cokernel of 
$ \Pic^0(X)\ra\Pic^0(\foX)$ is a torsion group and if
$X $ is numerically 
$\QQ $-factorial with respect to 
$S $, then 
$Y $ is numerically 
$\QQ $-factorial with respect to 
$T $.
\end{proposition}

\proof
Given a cycle 
$D\in Z^1(Y) $; then there is an integer 
$n>0 $ and a numerically trivial
Weil divisor 
$N  $ on $X$ such that 
$f^*(nD)+N $ is  Cartier   near 
$S $. Passing to a multiple and adding a numerically trivial
Cartier divisor, we can assume that  the corresponding reflexive 
$\O_{X} $-module 
$\shL  $ is trivial on 
$\foX $. Hence 
$\shM=f_*(\shL) $ is invertible at  
$T $ and is represented by 
$nD+f_*(N) $. QED.

\medskip
Let 
$nR\subset X $ be the infinitesimal neighborhoods of the
contracted curve
$R $. Then the inverse system of groups 
$n\mapsto H^1(R,\O_{nR}) $ is eventually constant, and we obtain a
well-defined group scheme of finite type 
$\Pic^0_{\foX/K} = \Pic^0_{nR/K}$, provided 
$n $ is sufficiently large. We can use the group scheme 
$G $ defined by the exact sequence 
$$
\Pic^0_{X/K}  \lra \Pic^0_{\foX/K}   \lra G  \lra 0
$$
to check the hypothesis of the previous result,  using the
following observation:

\begin{lemma}
\label{torsion cokernal}
Let 
$R\subset S $ be a closed subscheme of a proper 
$K $-scheme 
$S $. If 
$H^1(S,\O_{S})\ra H^1(R,\O_{R}) $ is surjective, or if the base
field
$ K$ is of characteristic 
$p>0 $ and the cokernel 
$G $ of the homomorphism
$\Pic^0_{S/K}  \ra \Pic^0_{R/K} $ is a  unipotent group scheme, 
then the cokernel of 
$ \Pic^0(S)\ra\Pic^0(R) $ is a torsion group.
\end{lemma}

\proof
Let 
$\shL $ be an numerically trivial invertible 
$\O_R $-module and 
$l\in \Pic^0_{R/K}   $ the corresponding {\em rational} point. The
map 
$H^1(S,\O_{S})\ra H^1(R,\O_{R}) $ is the tangential map for the
homomorphism 
$ \Pic^0_{S/K}  \ra \Pic^0_{R/K}$; if the first condition holds,
this map is surjective, and we find a {\em closed} point 
$m\in \Pic^0_{S/K} $ mapping to 
$l $.
Assume that the second condition holds. 
From the definition of unipotent group schemes
\cite{SGA 3}, SGA 3, p.~534, it follows immediately that 
$G(K) $ is a torsion group. Replacing 
$\shL $ by a multiple we also find a closed point 
$ m\in \Pic^0_{S/K}$ mapping to 
$ l$.

If the base field is algebraically closed, the point 
$m $ is represented by a numerically trivial invertible 
$\O_{S} $-module 
$\shM $ restricting to 
$\shL $.
In general,  we find a finite field extension 
$K\subset K' $ and an invertible sheaf 
$\shM' $ on 
$X'=X\otimes K' $ restricting to 
$\shL'=\shL\otimes K' $ on 
$R'=R\otimes K' $. Since  
$p_*(\O_{S'})$ is a locally free 
$\O_{S} $-module, say of rank 
$n >0$, there is a commutative diagram of norm homomorphism
\cite{EGA II}, EGA II, 6.6.8 
$$
\begin{CD}
\Pic(S')     @>>>  \Pic(R')   \\
@VN VV     @VV N V\\
\Pic(S)    @>>>   \Pic(R). 
\end{CD}
$$
Now 
$\shL^{\otimes n} = N(\shL') $ extends to 
$\shM=N(\shM') $. It remains to check that 
$\shM $ is numerically trivial. Using base change, we can assume
that 
$S $ is an integral curve. Then 
$\shM' $ is represented by a Cartier divisor 
$D' $ whose cycle 
$\sum n_{x'}x' $ satisfies 
$\sum n_{x'}\dim \k(x') =0$. Obviously, the cycle   
$\cyc(f_*(D)) $ is of degree zero. QED.

\begin{corollary}
\label{factorial and irregularity}
Let 
$Y $ be a proper normal surface over a  finite ground field $K$ 
or with 
$H^2(Y,\O_{Y})=0 $. Then 
$Y$ is numerically 
$\QQ $-factorial, and in particular projective.
\end{corollary}

\proof
Let 
$ X\ra Y $ be a resolution of singularities, and 
$\foX\subset X $ the formal completion of the contracted curve
$R\subset X$. 
If 
$H^2(Y,\O_{Y})=0 $ holds, then the restriction
$H^1(X,\O_{X})\ra H^1(\foX,\O_{\foX}) $ is surjective,  and the
claim follows from the above lemma and
\ref{factoriality and birational}.
If the ground field 
$K $ is finite, then 
$\Pic^0_{\foX/K} $ contains only finitely many rational points,
and 
$\Pic^0(\foX) $ must be a finite group.
Consequently, 
$Y $ is even 
$\QQ $-factorial. QED.

\medskip
Now we turn to the following situation: Let 
$R\subset X $ be a   curve on a proper normal surface 
$X $, and 
$A\in Z^1(X)$ a Weil divisor which is already Cartier near 
$R $ and numerically trivial on 
$R $. We would like to know whether we can make the {\em linear}
class of 
$A $ trivial on 
$R $. Here the canonical class 
$K_X $ becomes useful:

\begin{proposition}
\label{canonical cycle}
Let 
$A\in Z^1(X) $ be Cartier near 
$R $ and numerically trivial on 
$R $. If the curve 
$mR\subset X $ is contained in the base scheme of 
$K_X+mR $, then there is an integer 
$n>0 $ and a numerically trivial Cartier divisor 
$N $ such that the linear class of 
$nA+N $ is trivial on the infinitesimal neighborhood
$mR $.
\end{proposition}

\proof
We have an exact sequence 
$$
H^1(X,\O_{X})  \lra H^1(R,\O_{mR})  \lra  
H^2(X,\O_{X}(-mR)) \lra H^2(X,\O_{X}),
$$
and the map on the right is dual to the canonical  map 
$H^0(X,\O_{X}(K_X))\ra H^0(X,\O_{X}(K_X+mR)) $. Recall that the
base scheme of 
$K_X+mR $ is the intersection of all effective curves 
$C\subset X $ linearly equivalent to 
this class, hence the map is surjective. We conclude that 
$H^1(X,\O_{X})  \ra H^1(R,\O_{mR}) $ is surjective, and the claim
follows from 
\ref{torsion cokernal}. QED.

\medskip
In positive characteristics 
the infinitesimal neighborhoods are irrelevant for our problem:

\begin{proposition}
\label{trivial in characteristic p}
Let 
$A\in Z^1(X) $ be Cartier near 
$R $ and   trivial on 
$R $. Assume that the base field
$ K$ is of characteristic 
$p>0 $. Then the linear class of 
$p^mA $ is trivial on 
$mR $ for all integers 
$m>0 $.
\end{proposition}

\proof
We make induction on the integer
$ m$. Set $\shI_m=\O_{X}(-mR)$. The exact sequence 
$$
0 \lra \shI_m/\shI_{m+1} \lra  \O_{(m+1)R}^\times \lra
\O_{mR}^\times 
\lra 1 
$$
yields an exact sequence 
$$
H^1(R,\shI_m/\shI_{m+1}) \lra \Pic((m+1)R) \lra \Pic(mR) \lra 0. 
$$
By induction, 
$p^{m }A $ is trivial on 
$ mR$. Hence the class of 
$p^{m }A$ on
$(m+1)R $ lies in the image of 
$H^1(R,\shI_m/\shI_{m+1}) $. Since this group is annihilated  by 
$p $, we conclude that 
$ p^{m+1 }A$ is trivial on 
$(m+1)R $. QED.

\subsection{Applications}

In this   section we will collect some sufficient conditions
for contractibility and discuss an example. We assume that 
$X $ is a proper normal surface, 
$R\subset X $ is a connected  negative definite curve, and 
denote by 
$\foX \subset X $ the corresponding formal completion.  Our most
general criterion is the following:

\begin{theorem}
\label{sufficient condition}
Assume that 
$X $ is numerically 
$\QQ $-factorial with respect to 
$R $. If  the curve
$mR\subset X $ is contained in the base scheme of the class
$K_X+mR $ for all integers 
$m>0 $, or if the base field
$K $ is of characteristic 
$p>0 $ and the cokernel of 
$\Pic^0_{X/K}\ra\Pic^0_{R/K} $ is a unipotent group scheme, then 
$R $ is contractible.
\end{theorem}

\proof
Let 
$A  $ be a Weil divisor which is a support function of
the  pseudo-effective cone 
$\NEP(X) $ with respect to the extremal subcone 
$P=\sum\RR_+ R_i $ generated by the irreducible components 
$R_i\subset R $. Since 
$ X$ is numerically 
$\QQ $-factorial with respect to 
$R $, we can assume that 
$A $ is Cartier  near 
$R $.
 Given an integer 
$m>0 $. If the first condition holds, we invoke
\ref{canonical cycle} and find an integer 
$n>0 $ and a numerically trivial divisor 
$N\in \Div^0(X) $ such that the linear class of 
$nA+N $ is trivial on 
$mR $. If the second conditions holds, we first find such
$n $ and 
$N $ such that 
$nA+N $ is trivial on 
$R $. Hence 
$p^m(nA+N) $ is trivial on 
$mR $ by 
\ref{trivial in characteristic p}. Now the claim follows from
theorem
\ref{complementary cycles}. QED.

\medskip
The following will be useful for the classification of proper 
normal surfaces of Kodaira dimension $\kappa(X)=-\infty$:

\begin{corollary}
\label{contractibles and canonical cycle}
If the linear class 
$K_X+mR $ is not effective for all integers 
$m>0 $, or if the base field
$K $ is of characteristic 
$p>0 $ and 
$K_X+R $ is not effective, then 
$R $ is contractible.
\end{corollary}

\proof
Since 
$ K_X+R$ is not effective, the same holds for 
$ K_X$, and 
$H^2(X,\O_{X}) $ must vanish. According to
\ref{factorial and irregularity}, the surface 
$X $ is numerically 
$\QQ $-factorial.
If the first condition holds, the base scheme of 
$K_X+mR $ is the whole surface 
$X $, thus contains 
$mR $. If the second condition holds, the map
$H^1(X,\O_{X})\ra H^1(R,\O_R) $ is surjective, hence  
$\Pic^0_{X/K}\ra \Pic^0_{R/K} $ is an epimorphism. 
The claim now follows from 
\ref{sufficient condition}. QED.

\medskip
We can generalize the criterion of Castelnuovo-Enriques from
regular to normal proper surfaces:

\begin{corollary}
\label{Castelnuovo criterion}
Assume that 
$X $ is numerically 
$\QQ $-factorial with respect to 
$R $. If 
$R\subset  X $ is irreducible with
$R\cdot K_X\leq 0 $, then 
$R $ is contractible.    
\end{corollary}
  
\proof
Assume that 
$mR $ is not in the fixed scheme of 
$K_X+mR $ for some integer 
$m>0 $. Choosing 
$m $ minimal, we can represent 
$K_X+mR $ by a curve 
$C\subset X $ not containing 
$R $. But 
$$
0\geq K_X\cdot R = C\cdot R - mR^2 >0 
$$
gives a contradiction. According to
\ref{sufficient condition}, 
$R $ is contractible.   QED. 

\begin{remark}
For log-terminal surfaces 
$X $ and 
$K_X\cdot R<0 $, this becomes a special case of the contraction
theorem from Mori theory 
(\cite{Kawamata; Matsuda; Matsuki 1987}, theorem 3.2.1).
\end{remark}

\medskip
The following is already contained in
\cite{Artin 1962}, theorem 2.9:

\begin{proposition}
\label{finite ground field}
Assume that the 
$K $-surface
$X $ is already defined over a subfield 
$K'\subset K $, such that 
$K'$ is a finite field and that 
$K'\subset K $ is the composition of    purely transcendental
and radical extensions. Then every negative definite
curve 
$R\subset X $ is contractible.  
\end{proposition}

\proof
Let 
$X' $ be the normal surface over $K'$ with 
$X=X'\otimes_{K'}K $. According to 
\cite{SGA 6}, X, 7.17.4,
the mapping 
$Z^1(X',\QQ)\ra Z^1(X,\QQ) $ is surjective up to linear
equivalence. Hence there is a negative definite curve 
$R'\subset X' $ with 
$R=R'\otimes_{K'}K $, and we can assume that 
$K $ is finite. Then 
$X $ is 
$\QQ $-factorial, and 
$\Pic^0(mR)\subset\Pic^0_{mR/K}(K) $ are finite groups. According
to
\ref{complementary cycles}, the curve 
$R $ is contractible. QED.

\begin{example}
The case of geometrically ruled surfaces is already instructive.
Let 
$C $ be a normal proper connected curve and 
$p:X\ra C $ a geometrically ruled surface. Each section 
$R\subset X $ determines an extension 
$$
0 \lra \O_{C}  \lra \shE  \lra  \shL \lra 0 
$$
with 
$\shE=p_*(\O_{X}(R)) $ and 
$\shL=p_*(\O_R(R)) $. Hence we have 
$X=\PP(\shE) $ and 
$R=\PP(\shL) $, and  
$\O_R(R)\simeq p^*(\shL)\mid R$ holds. Assume there is a section
with
$R^2<0 $. Let 
$A $ be a divisor representing 
$p^*(\shL)\otimes\O_{X}(-R) $. Then 
$A $ is trivial on 
$R $, and in  characteristic 
$p>0 $ we deduce  from
\ref{complementary cycles} and 
\ref{trivial in characteristic p} that 
$R $ is contractible. However,  the situation is more
complicated in characteristic zero. We have 
$$
\Pic(2R)=\Pic(R)\oplus H^1(C,\shL^\vee); 
$$
decomposing the class of 
$A $ into 
$(0,\alpha) $, one can show that 
$\alpha $  corresponds to  the  Yoneda class in 
$\Ext^1(\shL,\O_{X}) $ of the extension 
$\shE $.  In characteristic zero, we conclude that 
$R $ is contractible if and only if the extension 
$\shE $ is split. In this case a splitting 
$\shE\ra \O_{C} $ yields another section 
$A\subset X $ disjoint to 
$R $ with 
$A^2>0 $, defining the contraction.
\end{example} 

\subsection{A finiteness result for models}

In this section we show that the {\itshape model}   of a normal
surface
$X $ defined by an arbitrary Weil divisor 
$D $ is a scheme of finite type. This generalizes results of 
Zariski \cite{Zariski 1962}, Proposition 11.5,
Fujita
\cite{Fujita 1981}, p.\ 235,  and Russo
\cite{Russo}.

\begin{atom}
\label{homogeneous spectrum}
Let 
$X $ be a proper normal surface and 
$D\in Z^1(X) $ a Weil divisor. This yields a graded 
$K $-algebra 
$$
R(X,D)= \oplus_{n\geq 0} H^0(X,nD), 
$$
which in turn defines a homogeneous spectrum
$$
P(X,D)= \Proj (R(X,D)). 
$$
According to 
\cite{EGA II}, EGA II, 3.7.4 the open subset 
$U\subset X $ of all points 
$x\in X $ such that there is a homogeneous 
$s\in R_+(X,D) $ with 
$s(x)\neq 0 $ is the largest open subset on which the canonical 
map 
$\O_{X}\otimes R(X,D)\ra \oplus_{n\geq 0} \O_{X}(nD) $ defines a
morphism 
$r:U\ra P(X,D) $. We call the scheme
$P(X,D) $ the 
{\itshape $D $-model } of the surface 
$X $. 

Using 
\cite{EGA I}, EGA I, 6.8.2, we see that for each homogeneous
element
$s\in R_+(X,D) $ the open subset 
$D_+(s)\subset P(X,D) $ equals the affine hull 
$X_s^\aff =\Spec\Gamma(X_s,\O_{X}) $, hence if 
$X_{s_1}\cup\ldots\cup X_{s_n} $ is a covering of 
$U $, then 
$D_+(s_1)\cup\ldots\cup D_+(s_n) $ is a covering of 
$r(U)\subset P(X,D) $. 
In this section we will prove the following 
\end{atom}

\begin{theorem}
\label{finiteness of models}
Let 
$X $ be a proper normal surface, and 
$D\in Z^1(X) $ an arbitrary Weil divisor. Then the model
$P(X,D) $ is a separated normal $K$-scheme   of finite type of
dimension 
$\leq 2 $.
\end{theorem}

It is well known that the algebra 
$R(X,D) $ might fail to be finitely generated,
see for example \cite{Zariski 1962}, p.\ 562. We will explain
the geometric reason for this below. 
First, we record the following

\begin{corollary}
\label{hilberts problem}
Let 
$U $ be a  normal surface, not necessarily proper. Then 
$\Gamma(U,\O_U) $ is an integrally closed
$K $-algebra of finite type and of dimension 
$\leq 2 $.
\end{corollary}

\proof
By the Nagata compactification theorem
\cite{Nagata 1962}, we can find a proper
scheme 
$X $ containing 
$U $ as an   open   subset. Making a blow-up and a
normalization we can assume that 
$X $ is a proper normal surface and that 
$D=X\setminus U $ is a  Cartier divisor. Now 
$U^\aff $ is an affine open subset of 
$P(X,D) $, and the claim follows from 
\ref{finiteness of models}. QED.

\begin{atom}
\label{preliminary assumptions}
\textit{Proof of theorem 
\ref{finiteness of models}.}
We start with some preliminary reductions. The 
$D$-model does not change if we replace 
$D $ by a positive multiple. 
If 
$ H^0(X,nD)=0$ holds for all integers 
$n>0 $, then 
$P(X,D) $ is empty, hence it suffices to treat the case that 
$D $ is effective.  For each 
$n>0 $ let 
$F_n\subset X $ be the fixed curve of 
$nD $. The effective Weil divisor 
$M_n=nD-F_n $ has no fixed curve. Following Kawamata
\cite{Kawamata 1997}, definition 1.1, we call such Weil divisors
\emph{movable}, and refer to 
$M_n$ as the movable part  of 
$nD $.
The canonical map 
$F_n:H^0(X,M_n)\ra H^0(X, nD) $ is bijective, and 
$M_n\cdot C\geq 0 $ holds for all curves 
$C\subset X $. Let 
$B_n\subset X $ be the base scheme of the class 
$nD $, that is the scheme-theoretical intersection of all curves 
$C\subset X $ representing the class
$nD $. Clearly, 
$B_n $ contains 
$F_n $ and all closed points 
$x\in X $  where 
$nD $ is not Cartier. Passing to a multiple of 
$D $, we can assume that
the supports of 
$B_n $ and 
$F_n $ are independent of 
$n $.

If 
$M_n=0 $ hold for all integer 
$n>0 $, we have 
$P(X,D)=X^\aff $, hence we can assume that 
$D $ has a non-zero movable part. 

Now let 
$F_n'\subset X $ be the union of all connected components 
$C\subset F_n $ with
$C\cdot M_n>0 $, and let
$F_n''= F_n - F_n'$ be the union of the remaining connected
components. Since 
$M_n $ is movable, the condition 
$M_n\cdot F_n=0 $ is equivalent to 
$F_n'=0 $. For simplicity, we   set
$M=M_1 $ and
$F=F_1 $.  The proof of 
\ref{finiteness of models} will be completed by the next two
propositions:
\end{atom}

\begin{proposition}
\label{1-dimensional models}
Under the above assumptions, the scheme 
$P(X,D)  $ is a normal projective curve if and only if 
$M^2=0 $ and 
$M\cdot F=0 $ holds.
\end{proposition}

\proof
Assume that 
$M^2=0 $ holds. Then for each curve 
$\sum\lambda_iC_i \subset X$ representing 
$M $ we must have 
$ M\cdot C_i=0$. Since 
$M $ is movable, it  can also be represented by a curve 
$C\subset X $ disjoint to 
$C_i $, and we deduce that 
$M $ must be a globally generated Cartier divisor. Thus 
$Y=P(X,M) $ is a proper curve, defining a fibration 
$f:X\ra Y $, and 
$M $ is the preimage of an ample class on 
$Y $. 

Now assume  that $M\cdot F =0 $ also holds; since 
$M $ is movable, 
$F' $ must vanish. Thus
$F=F'' $ is supported by certain fibres 
$X_y $ of 
$f:X\ra Y $. Assume that 
$\Supp(F_y)=\Supp(  X_y) $ holds for some point 
$y\in Y $. Decompose
$X_y=\sum\lambda_iE_i $ and 
$F_y=\sum \mu_i E_i $ into prime cycles.
Rearranging the terms, we can assume that 
$\lambda_1/\mu_1 \geq \lambda_i/\mu_i $
holds for all indices 
$i $. We have 
$$
\lambda_1 F_y = \mu_1 X_y + 
\sum_{i>1} (\lambda_1\mu_i - \mu_1\lambda_i)E_i,
$$
where the first summand 
$\mu_1 X_y$ is movable, whereas the second  summand
$\sum_{i>1} (\lambda_1\mu_i - \mu_1\lambda_i)E_i $ is effective. 
Hence 
$E_1 $ is not contained in the base locus of 
$\lambda_1 D $, but by our assumption 
\ref{preliminary assumptions}, the supports of 
$F_n $ are constant, contradiction. So 
$\Supp(F) $ contains no fibre of 
$f:X\ra Y $.
Consider the open subset
$V=X\setminus F $. I claim that the canonical map 
$\O_Y\ra f_*(\O_V) $ is bijective. This is local in 
$Y $; passing to the henselization   
$\O_{Y,y} \subset \O_{Y,y}^\sim $, we can contract the components 
$F_y\subset X_y $, according to 
\cite{Bosch-Luetkebohmert-Raynaud 1990},
Proposition 4, p.\ 169. By Serre's condition 
$(S_2) $ we have 
$f_*(\O_{V})\otimes\O_{Y,y}^\sim = \O_{Y,y}^\sim $.
Consequently, the open subsets 
$D_+(s)=X_s^\aff $  form a
covering of 
$Y $ when 
$s $ ranges over the homogeneous elements of 
$R_+(X,D) $, and we infer 
$ P(X,D)=Y$.

Conversely, assume that 
$Y= P(X,D) $ is a projective curve, and consider the morphism 
$r:U\ra Y $, where 
$U\subset X $ is the maximal open subset as in 
\ref{homogeneous spectrum}. If 
$M^2>0 $ holds, then 
$M\cup F' $ supports a curve ample on itself, and 
its complement  is almost affine
\ref{almost affine open subset}; thus 
$P(X,D) $ would contain a 2-dimensional open subset, 
contradiction. So 
$M^2=0 $ holds, and we obtain a fibration 
$f:X\ra Y $ extending $r$. If 
$F'\neq 0 $ holds, the generic fibre 
$X_\eta $   would be affine, contradicting 
\cite{EGA I}, EGA I, 9.3.4. Hence we have 
$M\cdot F=0 $. QED.

\begin{proposition}
\label{2-dimensional models}
With the assumptions in 
\ref{preliminary assumptions}, the model 
$P(X,D) $ is a normal surface if and only if either
$M^2>0 $, or 
$M^2=0 $ and 
$M\cdot F\neq 0 $ holds.
\end{proposition}

\proof
The condition is sufficient: Let 
$C\subset X $ be a curve  representing the class 
$M+F' $  and 
$V=X\setminus C $.
In case 
$M^2>0 $ the curve 
$C $ supports a curve ample on itself. If 
$M^2=0 $ holds, 
$Y=P(X,M) $ is a proper curve, and no connected component of
$F' $ is vertical with respect to the corresponding fibration 
$f:X\ra Y $. Again 
$C $ supports a curve ample on itself. In both cases 
$V $ is almost affine and 
$F''\subset V $ is contracted in   
$V^\aff $, thus $V^\aff $ is an open subset of 
$P(X,D) $ which is 2-dimensional and of finite type over 
$K $.
According to 
\ref{1-dimensional models}, we have 
$M_n^2>0 $ or $M_n^2=0$ and $M_n\cdot F_n\neq 0$ for all 
$n>0 $, consequently 
$P(X,D) $ is a 2-dimensional scheme of finite type.

Conversely, the condition is necessary by 
\ref{1-dimensional models}.   QED.

\medskip
If the scheme 
$P(X,D) $ is not proper, the algebra 
$R(X,D) $ is not finitely generated. Concerning this, we have the
following characterization:

\begin{proposition}
\label{proper models}
Assume that the model
$P(X,D) $ is a surface. With the assumptions in 
\ref{preliminary assumptions}, the surface
$P(X,D) $ is  proper if and only if 
$M $ is 
$\QQ $-Cartier and 
$M\cdot F=0 $ holds.
\end{proposition}

\proof
  Assume that 
$M $ is 
$\QQ $-Cartier and 
$M\cdot F=0 $  holds. Thus 
$F'=0 $ and $M^2>0$ holds, and 
$F''\subset X $ is a negative definite curve. Let 
$f:X\ra Y $ be the contraction of the negative definite curve 
$R\subset X $ orthogonal to 
$M $. Since 
$M $ is 
$\QQ $-Cartier, 
$Y=P(X,M) $ holds. Choose homogeneous
$s_1,\ldots,s_n \in R_+(X,M) $ with 
$X=X_{s_1}\cup\ldots\cup X_{s_n} $ and 
$R\subset X_{s_i} $. Thus we have 
$U_{s_i}^\aff=X_{s_i}^\aff $ and deduce 
$P(X,D)=P(X,M) $. 

Conversely, assume that 
$Y=P(X,D) $ is a  proper surface. Then each curve 
$C\subset X $  representing the class 
$M_n+F_n' $ is the support of a curve ample on itself, hence
$V=X\setminus C $  is almost affine, and 
$Y $ is covered by the open subsets 
$V^\aff $. Hence we can extend 
$r:U\ra Y $ to a proper morphism 
$r:V\ra Y $ on the larger open subset 
$V=U\cup F'' $, thus 
$X=V$ and 
$F'=0 $ hold, and 
$M\cdot F=0 $ follows. 
  Assume that 
$M $ is not 
$\QQ $-Cartier at some point 
$x\in X $; then 
$D $ is also not 
$\QQ $-Cartier at 
$x $, and we have
$x\not\in U\cup F''$, contradiction. QED.


\vspace{3em}
\noindent
        Mathematisches Institut\\
        Ruhr-Universit\"at\\
        44780 Bochum\\
        Germany\\ 
        E-mail: s.schroeer@ruhr-uni-bochum.de

\end{document}